\documentclass[12pt,a4paper]{amsart}
\usepackage{amsfonts}
\usepackage{epsfig}
\usepackage{graphicx}
\usepackage{amsmath}
\usepackage{amssymb}

\newcommand{\Rr}{\mathbb{R}}

\DeclareMathAlphabet{\mathpzc}{OT1}{pzc}{m}{it}

\newcounter{main}

\newtheorem{theorem}{Theorem}[section]

\newtheorem{corollary}[theorem]{Corollary}

\newtheorem{maintheorem}{Theorem}

\newcommand{\blanksquare}{\,\,\,$\sqcup\!\!\!\!\sqcap$}

\newcounter{example}
{{\stepcounter{example}}{\flushleft {\bf Example \arabic{example}:}}}%
{\par}

\title[A remark on the topological stability of symplectomorphisms]{A remark on the topological stability of symplectomorphisms}

\author[M. Bessa]{M\'{a}rio Bessa}

\address{Departamento de Matem\'atica da Universidade do Porto, 
Rua do Campo Alegre, 687, 
4169-007 Porto, Portugal and ESTGOH-Instituto Polit\'ecnico de Coimbra, Rua General Santos Costa, 3400-124 Oliveira do Hospital, Portugal}

\author[J. Rocha]{Jorge Rocha}
\address{Departamento de Matem\'atica da Universidade do Porto, 
Rua do Campo Alegre, 687, 
4169-007 Porto, Portugal}
\email{bessa@fc.up.pt}
\email{jrocha@fc.up.pt}

\begin{document}

\begin{abstract}
We prove that the $C^1$ interior of the set of all topologically stable $C^1$ symplectomorphisms is contained in the set of Anosov symplectomorphisms.
\end{abstract}

\maketitle

\section{Introduction: basic definitions and statement of the results}

Let $M$ be a $2d$-dimensional smooth manifold endowed with a symplectic form, say a closed and nondegenerate 2-form $\omega$. The pair $(M,\omega)$ is called a symplectic manifold which is also a volume manifold by Liouville's theorem. Let $\mu$ be the so-called Lebesgue measure associated to the volume form wedging $\omega$ $d$-times, i.e., $\omega^d=\omega\land\dots\land\omega$. On $M$ we also fix a Riemannian structure which induces a norm $\|\cdot\|$ on the tangent bundle $T_xM$. We will use the canonical norm of a bounded linear map $A$ given by $\|A\|=\sup_{\|v\|=1}\|A\,v\|$. By Darboux theorem (see e.g.~\cite[Theorem 1.18]{MZ}) there exists an atlas $\{h_j\colon U_j\to\Rr^{2d}\}$ satisfying $h_j^*\omega_0=\omega$ with 
\begin{equation}\label{canonical symplectic form}
\omega_0=\sum_{i=1}^d dy_i\land dy_{d+i}.
\end{equation}

A diffeomorphism $f\colon(M,\omega)\to(M,\omega)$  is called a \emph{symplectomorphism} if $f^*\omega=\omega$. Observe that, since $f^*\omega^d=\omega^d$, a symplectomorphism $f\colon M\to M$ preserves the Lebesgue measure $\mu$.

Symplectomorphisms arise naturally in the classical and analytical mechanics formalism as first return maps of hamiltonian vector
fields. Thus, it has long been one of the most interesting research fields in mathematical physics. We suggest the reference \cite{MZ} for more  details on general hamiltonian and symplectic theories. 

Let $(\text{Symp}_{\omega}^1(M), C^1)$ denote the set of all symplectomorphisms of class $C^1$  defined on $M$, topologized with the usual $C^1$ Whitney topology.

We say that $f\in\text{Symp}_{\omega}^1(M)$ is \emph{Anosov} if, there exist $\lambda\in(0,1)$ and $C>0$ such that the tangent vector bundle over $M$ splits into two $Df$-invariant subbundles $TM=E^u\oplus E^s$, with $\|Df^n|_{E^s}\|\leq C\lambda^n$ and $\|Df^{-n}|_{E^u}\|\leq C\lambda^n$. We observe that there are plenty Anosov diffeomorphisms which are not symplectic.

We say that $f \in \mathcal{F}_{\omega}^1(M)$ if there exists a neighborhood $\mathcal{V}$ of $f$ in $\text{Symp}_{\omega}^1(M)$ such that any $g\in\mathcal{V}$ has all the periodic orbits hyperbolic. We will use the following weighty result, which is a direct consequence of a theorem of Newhouse (see ~\cite[Theorem 1.1]{Ne}).

\begin{theorem}\label{AC}
If $f\in \mathcal{F}^1_{\omega}(M)$ then $f$ is Anosov. 
\end{theorem}
 
Let us recall that a periodic point $p$ of period $\pi$ is said to be \emph{1-elliptic} if the tangent map $Df^\pi(p)$ has two (non-real) norm one eigenvalues and the other eigenvalues have norm different from one.
Actually, unfolding homoclinic tangencies of symplectomorphisms $C^1$-far from the Anosov ones, Newhouse was able to prove that, $C^1$-generically, symplectomorphisms are either Anosov or else the $1$-elliptic periodic points of $f$ are dense in the whole manifold. We observe that the existence of $1$-elliptic periodic points is a sufficient condition to guarantee that the system is not structurally stable.

Given $f,g\in\text{Symp}_{\omega}^1(M)$ we say that $g$ is \emph{semiconjugated} to $f$ if there exists a continuous and onto map $h\colon M\rightarrow M$ such that, for all $x\in M$, we have the following conjugacy relation $h(g(x))=f(h(x))$. 

We say that $f$ is \emph{topologically stable} in $\text{Symp}_{\omega}^1(M)$ if, for any $\epsilon>0$, there exists $\delta>0$ such that for any $g\in \text{Symp}_{\omega}^1(M)$ $\delta$-$C^0$-close to $f$, there exists a semiconjugacy from $g$ to $f$, i.e., there exists $h\colon M\rightarrow M$ satisfying $h(g(x))=f(h(x))$ and $d(h(x),x)<\epsilon$, for all $x\in M$. Once again we emphasize that our definition of topological stability is restricted to the symplectomorphism setting and not to the broader space of volume-preserving (or even dissipative) diffeomorphisms. Let us denote by $\mathpzc{TS}_\omega(M)$ the subset of $\text{Symp}_{\omega}^1(M)$ formed by the topologically  stable symplectomorphisms.

The notion of topological stability was first introduced by Walters in (\cite{W}) proving that Anosov diffeomorphisms are topologically stable. Afterwards, in (\cite{Ni}), Nitecki proved that topological stability was a necessary condition to get Axiom A plus strong transversality. Later, in (\cite{Ro2}), Robinson proved that Morse-Smale flows are topologically stable. We point out that Hurley obtained necessary conditions for topological stability (see \cite{H,H2,H3}). In the early nineties Moriyasu (\cite{M}) proved that the $C^1$-interior of the set of all topologically stable diffeomorphisms is characterized as the set of all $C^1$-structurally stable diffeomorphisms. A few years ago Moriyasu, Sakai and Sumi (see~\cite{MSS}) proved that, if $X^t$ is a flow in the $C^1$ interior of the set of topologically stable flows, then $X^t$ satisfies the Axiom A and the strong transversality condition. Recently, in \cite{BeRo}, the authors proved a version of ~\cite{MSS} for the class of incompressible flows and also for volume-preserving diffeomorphisms.

The result in this note is a generalization of the theorems in \cite{MSS,BeRo} for symplectomorphisms. 

Given a set $A\subset \text{Symp}_{\omega}^1(M)$ let $int_{C^1}(A)$ denote the interior of $A$ in $\text{Symp}_{\omega}^1(M)$ with respect to the $C^1$-topology.

\begin{maintheorem}\label{teo1}
If $f\in int_{C^1}(\mathpzc{TS}_\omega(M))$ then $f$ is Anosov.
\end{maintheorem}

It is well known that Anosov diffeomorphisms impose severe topological restrictions to the manifold where they are supported. In fact the known examples of Anosov diffeomorphisms are supported in infranilmanifolds. An old conjecture of Smale states that any Anosov diffeomorphism is conjugated to an Anosov automorphism defined on an infranilmanifold. We end the introduction with this simple consequence of Theorem~\ref{teo1}.

\begin{corollary}
If $M$ does not support an Anosov diffeomorphism, then $$int_{C^1}(\mathpzc{TS}_\omega(M))=\emptyset.$$
\end{corollary}

\bigskip

\section{Proof of Theorem~\ref{teo1}}\label{diff}

Given $f\in int_{C^1}(\mathpzc{TS}_\omega(M))$, we prove that all its periodic orbits are hyperbolic; from this it follows that $int_{C^1}(\mathpzc{TS}_\omega(M))\subset \mathcal{F}^1_{\omega}(M)$. Then, using Theorem~\ref{AC}, we obtain that any $f\in int_{C^1}(\mathpzc{TS}_\omega(M))$ is an Anosov symplectomorphism. 

By contradiction let us assume that there exists $f\in int_{C^1}(\mathpzc{TS}_\omega(M))$ having a  non-hyperbolic orbit $p$ of period $\pi$.

In order to proceed with the argument we need to $C^1$-approximate the symplectomorphism $f$ by a new one, $f_1$, which, in the local coordinates mentioned in (\ref{canonical symplectic form}) and given by Darboux's theorem, is \emph{linear} in a neighborhood of the periodic orbit $p$. To perform this task, in the sympletic context, it is required more differentiability of the symplectomorphism (cf. ~\cite[Lemma 3.9]{AM}).

Therefore, if $f$ is of class $C^\infty$, take $g=f$, otherwise we use \cite{Z} in order to obtain a $C^\infty$ symplectomorphism $h\in int_{C^1}(\mathpzc{TS}_\mu(M))$, arbitrarily $C^1$-close to $f$, and such that $h$ has a periodic orbit $q$, close to $p$, with period $\pi$. We observe that $q$ may not be the analytic continuation of $p$ and this is precisely the case when $1$ is an eigenvalue of $Df^{\pi}(p)$. 

If $q$ is not hyperbolic take $g=h$. If $q$ is hyperbolic for $Dh^{\pi}(q)$, then, since $h$ is $C^1$-arbitrarily close to $f$, the distance between the spectrum of $Dh^{\pi}(q)$ and the unitary circle $\textbf{S}^1$ can be taken arbitrarily close to zero. This means that we are in the presence of weak hyperbolicity, thus in a position to apply ~\cite[Lemma 5.1]{HT} to obtain a new symplectomorphism $g\in\text{Symp}_{\omega}^\infty(M)\cap int_{C^1}(\mathpzc{TS}_\omega(M))$, $C^1$-close to $h$ and such that $q$ is  a non-hyperbolic periodic orbit.

Now, we use the \emph{weak pasting lemma} (\cite[Lemma 3.9]{AM})  in order to obtain $f_1 \in int_{C^1}(\mathpzc{TS}_\omega(M))$ such that, in local canonical coordinates, $f_1$ is linear and equal to $Dg$ in a neighborhood of the periodic non-hyperbolic orbit, $q$. Moreover, the existence of an eigenvalue, $\lambda$, with modulus equal to $1$ is associated to a sympletic invariant two-dimensional subspace. Furthermore, up to a perturbation using again ~\cite[Lemma 5.1]{HT}, $\lambda$ can be taken rational. This fact assures the existence of periodic orbits arbitrarily close to the $f_1$-orbit of $q$.

Finally, we $C^1$-approximate $f_1$ by $f_2\in int_{C^1}(\mathpzc{TS}_\omega(M))$ such that $q$ is a hyperbolic periodic orbit for $f_2$. This is a contradiction because $f_2$ is semiconjugated to $f_1$, although there is an $f_1$-periodic orbit (different from $q$) contained in a small neighborhood of the orbit of $q$ and the same cannot occur for $f_2$ because $q$ is a hyperbolic periodic orbit for $f_2$.

\section*{Acknowledgements}

The authors were partially supported by the FCT-Funda\c{c}\~ao para a Ci\^encia e a Tecnologia, project PTDC/MAT/099493/2008. MB was partially supported by FCT (SFRH/BPD/20890/2004).

\end{document}